\documentclass[a4paper, 10pt]{article}
\usepackage[latin1]{inputenc}
\usepackage{amscd, amssymb}
\usepackage[mathscr]{eucal}
\usepackage{mathrsfs} 
\usepackage{amsfonts}
\usepackage{amsmath}
\usepackage{amsthm}
\usepackage{latexsym}
\newcommand{\C}{\mathbf{C}}
\newcommand{\s}{\sigma}
\newcommand{\N}{\mathbf{N}}
\newcommand{\CU}{\mathbf{C}[[U]]}
\newcommand{\Ct}{\mathbf{C}[[t_1,...,t_n]]}
\newcommand{\pf}{\textsl{Proof.}}
\newcommand{\A}{A\otimes A^o}
\newcommand{\epf}{\hspace{14,15 cm} $\Box$}
\newcommand{\Wr}{\C \Gamma \sharp \C\left\langle x,y\right\rangle/(xy-yx=\Lambda)}

\newcommand{\G}{\mathbf{\Gamma}}
\newcommand{\g}{\gamma}
\newcommand{\h}{H_{t,k,c}(\G_N)}
\newcommand{\hh}{H_{1,k,c}(\G_N)}
\newcommand{\End}{\mathrm{End}}
\newcommand{\Ext}{\mathrm{Ext}}
\newcommand{\Hom}{\mathrm{Hom}}
\newcommand{\Tr}{\mathrm{Tr}}
\newcommand{\im}{\mathrm{Im}}
\newcommand{\Ker}{\mathrm{Ker}}
\newcommand{\rk}{\mathrm{rk}}
\newcommand{\ind}{\mathrm{Ind}}
\newcommand{\di}{\mathrm{dim}\,}
\newcommand{\vn}{\vec{N}}
\newtheorem{defi}{Definition}[section]
\newtheorem{thm}[defi]{Theorem}
\newtheorem{lem}[defi]{Lemma}
\newtheorem{prop}[defi]{Proposition}
\newtheorem{coro}[defi]{Corollary}
\newcommand{\rem}{\textbf{Remark.}}
\newcommand{\Z}{\mathbf{Z}}
\usepackage{latexsym}
\begin{document}

\title{Finite dimensional representations of symplectic
reflection algebras associated to wreath products II}

\author{Silvia Montarani}

\maketitle

\section{Introduction and main result}\label{intro}
In this paper we extend some results of a previous paper (\cite{EM}) about finite   dimensional representations of
the wreath product symplectic reflection algebra
$H_{1,k,c}(\G_N)$ of rank $N$ attached to the group
$\G_N=S_N\ltimes \Gamma^N$ (\cite{EG}), where $\Gamma\subset
SL(2,\C)$ is a finite subgroup, and $(k,c)\in C(\mathcal{S})$, where
$C(\mathcal{S})$ is the space of (complex valued ) class
functions on the set $\mathcal{S}$ of symplectic reflections of
$\G_N$. 

As in the previous paper, we use the results obtained in \cite{CBH} for the rank $1$ case and the homological properties of $H_{1,0,c}(\G_N)$ to obtain irreducible representations for small values of  $k\neq 0$. 
Specifically, let's denote by $B:=H_{1,c}(\Gamma)$ the rank one symplectic reflection algebra attached to the pair ($\Gamma$, $c$) (in the rank $1$ case there is no parameter $k$). For $\vn=(N_1,N_2,\dots,N_r)$ a partition of $N$,  we consider   $W=W_1\otimes\dots\otimes W_r$, an irreducible representation of  $S_{\vn}:=S_{N_1}\times\dots\times S_{N_r}$ such that the Young diagram of $W_i$ is a rectangle for any $i$ and $\{Y_i\}$, $i=1,\dots,r$,  a collection of irreducible non-isomorphic representations of   $B$  with $Ext^1_B(Y_i,Y_j)=0$ for any $i\neq j$. We denote by  $Y$ the tensor product $Y_1^{\otimes N_1}\otimes\dots\otimes Y_r^{\otimes N_r}$. Then $M'=W\otimes Y$ is in a natural way an irreducible representation of the algebra $S_{\vn}\sharp B^{\otimes N}\subset S_N\sharp B^{\otimes N}=H_{1,0,c}(\G_N)$. We show that the induced representation  $M=\ind_{S_{\vn}\sharp B^{\otimes N}}^{S_N\sharp B^{\otimes N}} M'$ of $H_{1,0,c}(\G)$ can be uniquely deformed along a linear subspace of codimension $r$ in $C(\mathcal{S})$.

\subsection{Symplectic reflection algebras of wreath product type}\label{wreprod}~Before stating our main result in more detail, we recall the definition of the  wreath product symplectic reflection algebra $H_{1,k,c}(\G_N)$.
Let $L$
be a $2$-dimensional complex vector space with a symplectic form
$\omega_L$, and consider the space $V=L^{\oplus N}$,
endowed with the induced symplectic form
$\omega_V={\omega_L}^{\oplus N}$. Let $\Gamma$ be a finite
subgroup of $Sp(L)$, and let $S_N$ be the symmetric group acting
on $V$ by permuting the factors. The group
$\G_{N}:=S_N\ltimes\Gamma^N \subset Sp(V)$ acts naturally on
$V$. In the sequel we will write $\gamma_i\in \G_N$ for any
element $\gamma\in \Gamma$ seen as an element in the i-th factor
$\Gamma$ of $\G_N$. The symplectic reflections in $\G_N$ are the
elements $s$ such that $\rk(Id-s)|_V=2$. $\G_N$ acts by conjugation
on the set $\mathcal{S}$ of its symplectic reflections. It is
easy to see that there are symplectic reflections
of two types in $\G_N$:

\medskip

(S) the elements $s_{ij}\g_i {\g_j}^{-1}$ where $i,j\in [1,N]$, $s_{ij}$ is the transposition $(ij)\in S_N$, and $\g\in\Gamma$,

($\Gamma$) the elements $\g_i$, for $i\in[1,N]$ and $\g\in\Gamma\backslash\{1\}$.
\medskip

Elements of type (S) are all in the same conjugacy class, while
elements of type ($\Gamma$) form one conjugacy class for any
conjugacy class of $\g$ in $\Gamma$. 
Thus elements $f\in \C[\mathcal{S}]$ can be written as pairs
$(k,c)$, where $k$ is a number (the value of $f$ on elements of
type (S)), and $c$ is a conjugation invariant function
on $\Gamma\setminus \lbrace{1\rbrace}$ (encoding the values of
$f$ on elements of type ($\Gamma$)).
    
For any $s\in\mathcal{S}$ we write $\omega_s$ for
the bilinear form on $V$ that coincides with $\omega_V$ on
$\im(Id-s)$ and has $\Ker(Id-s)$ as radical. Denote by $TV$ the
tensor algebra of $V$.

\begin{defi} \label{dewre}
 For any $f=(k,c)\in \C[\mathcal{S}]$, the symplectic reflection algebra $H_{1,k,c}(\G_N)$ is  the quotient $$
(\G_N\sharp TV)/\left\langle [u,v]-\kappa(u,v)\right\rangle_{u,v\in V}
$$ 
where 
$$
\kappa:V\otimes V\longrightarrow\C[\G_N]: (u,v)\mapsto \omega(u,v)\cdot 1 +\sum_{s\in \mathcal{S}}f_s\cdot \omega_s(u,v)\cdot s
$$
 with $f_s=f(s)$, and $\left\langle \dots\right\rangle$ is the two-sided ideal in 
the smash product $\G_N\sharp TV$ generated by the elements
$[u,v]-\kappa(u,v)$ for $u,v\in V$.

\end{defi} 

Following \cite{GG}, we will now  give a more explicit definition  of the algebra $\hh$ by generators and relations that we will need in the sequel.
It is clear that choosing a symplectic basis $x$, $y$ for $L$ we
can consider $\Gamma$ as a subgroup of $SL(2,\C)$. We will denote by $x_i$, $y_i$ the corresponding vectors in the i-th $L$-factor of $V=L^{\oplus N}$ and by $c_{\gamma}$ the value of the function $c$ on an  element $\g\in \Gamma$.

\begin{lem}([GG]) \label{pres}
The algebra $\hh$ is the quotient
of $\G_N\sharp TV$ by the following relations:
\begin{itemize}
\item[{\emph{(R1)}}] 
For any $i\in [1,N]$: $$[x_{i}, y_{i}]
=  1+ \frac{k}{2} \sum_{j\neq i}\sum_{\g\in\Gamma}
s_{ij}\g_{i}\g_{j}^{-1} + \sum_{\g\in\Gamma\backslash \{1\}}
c_{\g}\g_{i}\,.$$
\item[{\emph{(R2)}}] 
For any $u,v\in L$ and $i\neq j$:
$$[u_{i},v_{j}]= -\frac{k}{2} \sum_{\g\in\Gamma} \omega_{L}(\g u,v)
s_{ij}\g_{i}\g_{j}^{-1} \,.$$
\end{itemize}
\end{lem} 

\epf

In the case $N=1$, there is no parameter $k$ (there are no symplectic reflections of type (S)) and  

\begin{equation}\label{rankone}
H_{1,c}(\Gamma)=\Wr 
\end{equation}

where
\begin{equation}\label{centralel}
\Lambda=\Lambda(c)=1+\sum_{\g\in\Gamma\backslash\{1\}}c_{\g}\g\in
Z(\C[\Gamma])
\end{equation} 
is the central element corresponding to $c$.

\subsection{The main result}
\label{mainteo} 

Let $\{Y_i\}$ be a collection of non-isomorphic  irreducible representations 
of the algebra $B=H_{1,c}(\Gamma)$ for some $c$
\footnote{Such representations exist only for special $c$, as for
generic $c$ the algebra $H_{1,c}(\Gamma)$ is simple; see \cite{CBH}.} and let $Y=Y_1^{\otimes N_1}\otimes\dots\otimes Y_r^{\otimes N_r}$. 
Let $W=W_1\otimes\dots\otimes W_r$ be an irreducible representation of $S_{\vn}=S_{N_1}\times\dots\times S_{N_r}\subset S_N$, where $S_{N_i}$ is the subgroup of permutations moving only the indices $\{\sum_{j=1}^{i-1}N_j+1,\dots ,\sum_{j=1}^iN_j\}$.
There is a natural action of the algebra $S_{\vn}\sharp B^{\otimes N}$ on the vector space $M':=W\otimes Y$.
Namely, each copy of $H_{1,c}(\Gamma)$ acts in the corresponding
 $Y_i$, while $S_{\vn}$ acts on $W$ and simultaneously permutes
the factors in the product $Y$. We will denote this
representation by $M'_c$. Since the algebra $H_{1,0,c}(\G_N)$ is naturally isomorphic to $S_N\sharp B^{\otimes N}\supset S_{\vn}\sharp B^{\otimes N} $ we can form the induced module $M_c=\ind_{ S_{\vn}\sharp B^{\otimes N}}^{S_N\sharp B^{\otimes N}}M_c'$. The main theorem tells us when such a representation can be deformed to nonzero values of $k$. 

If the Young diagram of $W_i$ is a rectangle 
of height $l_i$ and width $m_i=N/l_i$ (the trivial representation
corresponds to the horizontal strip of height $1$) let ${\mathcal H}_{Y_i,m_i,l_i}$ be the hyperplane in $C(\mathcal S)$ 
consisting of all pairs $(k,c)$ satisfying the equation 
\begin{equation}\label{hyp}
\di Y_i+\frac{k}{2}|\Gamma|(m_i-l_i)+\sum_{\gamma\in
\Gamma\setminus \lbrace{1\rbrace}}c_\gamma\chi_{Y_i}(\gamma)=0,
\end{equation}
where $\chi_{Y_i}$ is the character of $Y_i$. 

Let $X=X(Y,W)$ be the moduli space of irreducible representations
of $H_{1,k,c}(\G_N)$  isomorphic to $M$ 
as $\G_N$-modules (where $(k,c)$ are allowed to vary). 
This is a quasi-affine algebraic variety:
it is the quotient of the quasi-affine variety
$\widetilde{X}(Y,W)$ of extensions of 
the $\G_N$-module $M$ 
to an irreducible $H_{1,k,c}(\G_N)$-module 
by a free action of the reductive group $G$ of basis changes in $M$
compatible with $\G_N$ modulo scalars. Let $f: X\to
C({\mathcal S})$ 
be the morphism which sends a representation to the
corresponding values of $(k,c)$. 

The main result of this paper is the following theorem that implies Theorem $3.1$ of \cite{EM} as a particular case.

\begin{thm}\label{mainth} Suppose $W_i$ has rectangular Young diagram for any $i=1,\dots,r$, and  $\Ext_{B}^1(Y_i,Y_j)=0$  for any $i\neq j$. Then:

(i) For any $c_0$ the representation $M_{c_0}$ of $H_{1,0,c_0}( \G_N)$ 
can be formally deformed to a representation of $H_{1,k,c}(\G_N)$ along the codimension $r$ linear subspace $\displaystyle{\bigcap_{i=1}^r{\mathcal H}_{Y_i,m_i,l_i}}$, but not
in other directions. This deformation is unique. 

(ii) The morphism $f$ maps $X$ to $\displaystyle{\bigcap_{i=1}^r{\mathcal H}_{Y_i,m_i,l_i}}$  and is \'etale at $M_{c_0}$ for all $c_0$. Its restriction to the formal
neighborhood of $M_{c_0}$ is the deformation from (i).  

(iii) There exists a nonempty Zariski open subset 
$\mathcal U$ of $\displaystyle{\bigcap_{i=1}^r{\mathcal H}_{Y_i,m_i,l_i}}$ such that for $(k,c)\in
\mathcal U$, the algebra $H_{1,k,c}(\G_N)$ admits a 
finite dimensional irreducible representation isomorphic to $M$ as a 
$ \G_N$-module.  
\end{thm}

The proof of this theorem occupies the remaining sections of the paper.

\section{Main tools for the proof of Theorem \ref{mainth}} 
In this section we list some results we'll need to prove our main theorem.

\subsection{Representations of $S_N$ with rectangular Young diagram} We will need the following standard results from  the representation
theory of the symmetric group. The proofs are well known (for a quick reference see the  previous article \cite{EM}). Denote by $\mathfrak{h}$ the
reflection representation of $S_N$. For a Young diagram $\mu$ we
denote by $\pi_{\mu}$ the corresponding irreducible
representation of $S_N$. 

\begin{lem}\label{sym}
\begin{center}
\begin{itemize}
\item[{(i)}] $\Hom_{S_N}(\mathfrak{h}\otimes \pi_{\mu},\pi_{\mu})=\C^{m-1}$, where $m$ is the number of corners of the Young diagram $\mu$. In particular  $\Hom_{S_N}(\mathfrak{h}\otimes \pi_{\mu},\pi_{\mu})=0$ if and only if $\mu$ is a rectangle.
\item[{(ii)}] The element $C=s_{12}+s_{13}+ \cdots +s_{1n}$ acts by a scalar in $\pi_{\mu}$ if and only if $\mu$ is a rectangle.
 In this case the scalar corresponding to $C$ is the content  $\mathbf{c}(\mu)$ of the Young diagram $\mu$ i.e. the sum of the signed distances of the cells from the diagonal.
\end{itemize}
\end{center}
\end{lem}

\epf

\subsection{Irreducible representations of $H_{1,c}(\Gamma)$}\label{rank1}

We recall some  results of \cite{CBH} about the irreducible representations of the rank $1$ algebra $H_{1,c}(\Gamma)$. In \cite{CBH}, the classification of the finite-dimensional irreducible representations of   $H_{1,c}(\Gamma)$ is obtained by establishing a Morita equivalence of this algebra with a deformation of the path algebra of the double $\overline{Q}$ of the quiver $Q$, associated to the group $\Gamma$ in the McKay correspondence. This algebra is called \emph{deformed preprojective algebra}, and we will denote it by $\Pi^{\lambda}(Q)$. Here $\lambda\in \C^{\nu}$, where $\nu$ is the number of non isomorphic irreducible representations of $\Gamma$, is a parameter related to our parameter $c$ in the following way. If we denote by $\{V_i\}$, $i=1,\dots,\nu$ the collection of irreducible non isomorphic representations of $\Gamma$, with $V_1$ corresponding to the trivial representation, and by $\Lambda=1+\sum_{\g\in\Gamma\backslash\{1\}}c_{\g} \g$ the central element of $\C[\Gamma]$ that appears in the definition of $H_{1,c}(\Gamma)$ given at the end of Section \ref{wreprod} (see equation (\ref{centralel})), then we have:
$$\lambda_i=\Tr|_{V_i}\,\Lambda.$$ Thus, in particular, the  representations of $H_{1,c}(\Gamma)$ can be seen as representations of the path algebra of a certain quiver, so that we can attach a \emph{dimension vector} to each of them. We recall that one can associate an affine root system  to the McKay graph $Q$. The roots of such system can be distinguished into real roots and imaginary roots. The real roots are divided into positive and negative roots, and are the images of the coordinate vectors of $\Z^{\nu}$ under sequences of some suitably defined reflections. The  imaginary roots, instead, are all the non-zero integer multiples of the vector $\delta$, with $\delta_i=\di V_i$. If ''$\cdot$'' is  the  standard scalar product in $\C^{\nu}$, let's denote by $R_{\lambda}$ the set of real roots $\alpha$  such that $\lambda\cdot \alpha=0$, and by $\Sigma_{\lambda}$ the unique basis of $R_{\lambda}$ consisting of positive roots.
We have the following result:
\begin{thm}[\cite{CBH}, Theorem $7.4$]\label{cbht}
If $\lambda\cdot \delta\neq 0$, then $\Pi^{\lambda}(Q)$ has only finitely many finite-dimensional simple modules up to isomorphisms, and they are in one-to-one correspondence with the set $\Sigma_{\lambda}$. The correspondence is the one assigning to each module its dimension vector.
\end{thm}  

\epf

We recall that $\lambda\cdot \delta=\Tr|_{R} \Lambda$, where $R$ is the regular representation of $\Gamma$. Thus the theorem holds in the case $\Tr|_{R} \Lambda\neq 0$. Looking at our definition of $\Lambda$ we deduce that this is exactly our case.

\subsection{A standard proposition in deformation theory}

In this section we describe the cohomological strategy we will use to prove Theorem \ref{mainth}. Let $A$ be an associative algebra over $\C$. In what follows, for each $A$-bimodule $E$, we write $H^n(A,E)$ for  the $n$-th Hochschild cohomology group of $A$ with coefficients in $E$. 
We remark that    $H^i(A,E)$  coincides with the vector space $Ext^i_{\A}(A, E)$, where $A^o$ is the opposite algebra of $A$.

Let $A_U$ be a flat formal deformation of $A$ over the formal
neighborhood of zero in a
finite dimensional vector space $U$ with coordinates
$t_1,...,t_n$. This means that $A_U$ is an algebra over $\CU =
\Ct$ which is \textsl{topologically free} as a $\CU$-module (i.e., $A_U$  is isomorphic as a $\CU$-module to $A[[U]]$), together with a fixed isomorphism of algebras $A_U/JA_U\cong A$, where $J$ is the maximal ideal in $\CU$.  
Given such a deformation, we have a natural linear map\  $\phi:U \longrightarrow H^2(A, A)$.

Explicitly, we can think of $A_U$ as $A[[t_1,...,t_n]]$ equipped
with a new $\Ct$-linear (and continuous) associative product defined by: 
$$
a \ast b =\sum_{p_1,...,p_n}c_{p_1,...,p_n}(a,b)\prod_j t_j^{p_j} \qquad a,b \in A  
$$
where $c_{p_1,...,p_n}:A\times A\longrightarrow A$ are
$\C$-bilinear functions and $c_{0,...,0}(a,b)=ab$, for any
$a,b\in A$. 

Imposing the associativity condition on $\ast$, one obtains that $c_{0,...,1_j,...0}$ must be Hochschild $2$-cocycles for each $j$.  
The map $\phi$ is given by the assignment
$(t_1,\cdots,t_N)\rightarrow \sum_jt_j\,[c_{0,...,1_j,...0}]$ for
any $(t_1,\cdots,t_n)\in U$, where $[C]$ stands for the
cohomology class of a cocycle $C$. 

Now let $M$ be a representation of $A$. In general it does not
deform to a representation of $A_U$. However we have the
following standard proposition. Let $\eta: U\to H^2(A, \End M)$ be the composition of $\phi$ with the natural map $\psi:H^2(A,A)\longrightarrow H^2(A, \End\,M)$. 

\begin{prop}[\cite{EM}] \label{defo}
Assume that $\eta$ is surjective with kernel $K$, and $H^1(A,\End M)=0$. Then: 

(i) There
exists a unique smooth formal subscheme $S$ of the formal neighborhood
of the origin in 
$U$, with tangent space $K$ at the origin, 
such that  $M$ deforms to a representation of the algebra 
$A_S:=A_U \hat\otimes_{\CU}\C[S]$ (where $\hat\otimes$ is the
completed tensor product).

(ii) The deformation of $M$ over $S$ is unique.
\end{prop}
 
 \epf

We will show in the next sections that  the conditions of Proposition \ref{defo} hold in our case.

\subsection{Homological properties of the algebra $H_{1,c}(\Gamma)$.}

In order to show that we can make use of  Proposition \ref{defo} we will need  the following definitions and results.

\begin{defi}[ \cite{VB1,VB2,EO}]\label{dclass}
An algebra $A$ is defined to be in the class $VB(d)$ if it is of finite Hochschild dimension (i.e. there exists $n\in\N$ s.t. $H^i(A,E)=0$ for any $i>n$ and any $A$-bimodule $E$) and $H^{\ast}(A,A\otimes A^o)$ is concentrated in degree $d$, where it is isomorphic to $A$ as an $A$-bimodule.
\end{defi}
The nice properties of the class  $VB(d)$  are clarified by the following
result by Van den Bergh.

\begin{thm}[\cite{VB1,VB2}]\label{vdb}
If $A\in VB(d)$ then for any $A$-bimodule $E$, 
the Hochschild homology $H_i(A,E)$ is isomorphic to the
Hochschild cohomology $H^{d-i}(A,E)$. 
\end{thm}

\epf

\rem\ The isomorphism of Theorem \ref{vdb} depends on the choice of a bimodule isomorphism \\$\phi: H^d(A, A\otimes A^o)\stackrel{\sim}{\to} A\,$  and is unique for any such choice.
\vspace{0.4 cm}

Now let $B=H_{1,c}(\Gamma)$ and let $Y_i$ be  any irreducible representation of $B$. The following proposition  summarizes the homological properties of this algebra that we will use in  the proof of our main theorem. All the proofs can be found in the  previous article (\cite{EM}). 

\begin{prop}\label{vander}

\begin{enumerate}
\item[]\bigskip
\item[(i)] The algebra $B$ belongs to the class $VB(2)$;
\item[(ii)] $H^2(B,\End\,Y_i)=H_0(B, \End\, Y_i)=\C$;
\item[(iii)] $H^1(B,\End Y_i)=0$.
\end{enumerate}
\end{prop}

\epf

\section{Proof of Theorem \ref{mainth}}
\subsection{Homological properties of
$A=H_{1,0,c}(\G_N)$.} 

We are now ready to prove Theorem \ref{mainth}. Let $A$ denote the algebra $H_{1,0,c_0}(\G_N)$. 
The algebra $A$ has a flat deformation over 
$U=C({\mathcal S})$, which is given by the algebra 
$H_{1,k,c_0+c'}(\G_N)$. The fact that this deformation is
flat follows from the so called \emph{ Poincar\'e-Birkhoff-Witt Theorem}  for symplectic reflection algebras (see \cite{EG}, Theorem $1.3$). Our job is to compute the cohomology groups $H^2(A,\End\,M)$, $H^1(A,\End\,M)$ and the map $\eta:U\longrightarrow H^2(A,\End\,M)$.

\begin{prop}\label{endot}
If  $W=W_1\otimes\dots\otimes W_N$ is an irreducible representation of $S_{\vn}$ such that $W_i$ has rectangular Young diagram for any $i$ and $Y=Y_1^{\otimes N_1}\otimes\dots\otimes Y_r^{\otimes N_r}$ is an irreducible representation of $B^{\otimes N}$ with $Y_i\ncong Y_j$ and \ $\Ext^1_B(Y_i, Y_j)=\Ext^1_{B\otimes B^{o}}(B,\Hom(Y_i,Y_j))=0$ for any $i\neq j$ then 
$$
H^2(A,\End\,M)=\bigoplus_i H^2(B, \End\,Y_i)=\C^r.
$$
\end{prop}

\pf\  The second equality follows from Proposition \ref{vander}, $(ii)$.
Let us prove the first equality.
We have: 
$$
H^{\ast}(A,\End\,M)=\Ext^{\ast}_{A\otimes A^o}(A,\End\,M)=
$$
$$
=\Ext^{\ast}_{S_N\sharp B^{\otimes N}\otimes S_N\sharp {B^o}^{\otimes N}}(S_N\sharp B^{\otimes N},\End\,M)=
$$
$$
=\Ext^{\ast}_{S_N\times S_N\sharp(B^{\otimes N}\otimes {B^o}^{\otimes N})}(S_N\sharp B^{\otimes N},\End\, M).
$$
Now, the  $S_N\times S_N\sharp(B^{\otimes N}\otimes {B^o}^{\otimes
N})$-module $S_N\sharp B^{\otimes N}$ is induced from the module $B^{\otimes N}$ over the subalgebra $S_N\sharp \left(B^{\otimes N}\otimes {B^o}^{\otimes N}\right)$, in which $S_N$ acts simultaneously  permuting the factors of $B^{\otimes N}$ and ${B^o}^{\otimes N}$ (note that $S_N\sharp (B^{\otimes N}\otimes {B^o}^{\otimes N})$ is indeed a subalgebra of $S_N\times S_N\sharp(B^{\otimes N}\otimes {B^o}^{\otimes N})$ as it can be identified with the subalgebra $D\sharp(B^{\otimes N}\otimes {B^o}^{\otimes N})$ where $D=\left\{(\sigma,\sigma),\; \sigma\in S_N\right\}\subset S_N\times S_N$). Applying the Shapiro Lemma, we get:
$$
 \Ext^{\ast}_{S_N\times S_N\sharp(B^{\otimes N}\otimes {B^o}^{\otimes N})}(S_N\sharp B^{\otimes N},\End\,M)=
$$
$$
=\Ext^{\ast}_{S_N\sharp(B^{\otimes N}\otimes {B^o}^{\otimes N})}(B^{\otimes N},\End\,M)=
$$
$$
={\left(\Ext^{\ast}_{B^{\otimes N}\otimes {B^o}^{\otimes N}}(B^{\otimes N},\End\, M)\right)}^{S_N}.
$$

 We observe now that the subalgebra $B^{\otimes N}$ is stable under the inner automorphisms induced by the elements $\s\in S_N\subset A$. Thus the induced $A$-module $M$ can be written as:
\begin{equation}\label{deco}
M=\s_1 M'\oplus \s_2 M'\oplus\dots\oplus\s_n M'
\end{equation}
where $n=\frac{N!}{N_1!\dots N_r!}$ and $\{\s_1,\dots,\s_n\}$ is a set of representatives for the left cosets of $S_{\vn}$ in $S_N$. Each direct summand $\s_l M'$ can itself be written as:
$$
\s_l M'=\s_l W\otimes \s_l Y=\s_l W\otimes Y_{\alpha_{\s_l^{-1}(1)}}\otimes\dots\otimes Y_{\alpha_{\s_l^{-1}(N)}}.
$$
The action of an element $\s (b_1\otimes\dots \otimes b_N)$ on a vector $\left(w\otimes y_1\otimes\dots\otimes y_N\right)_l \in \s_l M'$ is the following:
\begin{equation}\nonumber
\s (b_1\otimes\dots \otimes b_N)\left(w\otimes y_1\otimes\dots\otimes y_N \right)_l=\end{equation} 
\begin{equation}\label{action}
=\left(\s'w\otimes b_{(\s'\s_h)^{-1}(1)} y_{(\s'\s_h)^{-1}(1)}\otimes\dots\otimes b_{(\s'\s_h)^{-1}(N)} y_{(\s'\s_h)^{-1}(N)}\right)_h\in \s_hM'
\end{equation}
where $\s'\in S_{\vn}$,  $\s_h\in \{\s_1,\dots, \s_n\}$ are  the only elements satisfying $\s\s_l=\s_h\s'$. 
Thus the latter module equals:
$$
\left(\Ext^{\ast}_{B^{\otimes N}\otimes {B^o}^{\otimes N}}\left(B^{\otimes N},\End\left(\bigoplus_{l}\s_l M'\right)\right)\right)^{S_N}=
$$
$$
=\left(\Ext^{\ast}_{B^{\otimes N}\otimes {B^o}^{\otimes N}}\left(B^{\otimes N},\bigoplus_{l,h} \Hom(\s_l M', \s_h M')\right)\right)^{S_N}
$$
Since  the action of $B^{\otimes N}\otimes {B^o}^{\otimes N}$ does not permute the direct factors in $\displaystyle{\bigoplus_l \s_l M'}$ and is trivial on $W$ in the module $M'$ we can rewrite the last term as: 
$$
\left(\bigoplus_{l,h}\Ext^{\ast}_{B^{\otimes N}\otimes {B^o}^{\otimes N}}\left(B^{\otimes N},\Hom(\s_l M',\s_h M')\right)\right)^{S_N}=
$$
$$
\left(\bigoplus_{l,h}\Ext^{\ast}_{B^{\otimes N}\otimes {B^o}^{\otimes N}}\left(B^{\otimes N},\Hom(\s_l W,\s_h W)\otimes\Hom(\s_l Y,\s_h Y) \right)\right)^{S_N}=$$
$$
=\left(\bigoplus_{l,h}\Hom(\s_l W,\s_h W)\otimes \Ext^{\ast}_{B^{\otimes N}\otimes {B^o}^{\otimes N}}\left(B^{\otimes N},\Hom(\s_l Y,\s_h Y) \right)\right)^{S_N}=
$$
$$
=\left(\bigoplus_{l,h}\Hom(\s_l W,\s_h W)\otimes \Ext^{\ast}_{B^{\otimes N}\otimes {B^o}^{\otimes N}}\left(B^{\otimes N},\bigotimes_{i=1}^N\Hom( Y_{\alpha_{\s_l^{-1}(i)}},Y_{\alpha_{\s_h^{-1}(i)}})\right)\right)^{S_N}.
$$
We want now to  apply   the K\"unneth formula in degree $2$. Proposition \ref{vander}, $(iii)$ and our conditions on the $Y_i$s guarantee $Ext^1_B(Y_i,Y_j)=0$. Moreover, since for any $i\neq j$  $Y_i$, $Y_j$ are non-isomorphic irreducible representations of $B$ we have $\Ext^0_B(Y_i,Y_j)=\Hom_B(Y_i,Y_j)=0$.
We get:
$$
\left(\bigoplus_{\begin{array}{c} l,h\\
(\s_l^{-1}(i) \s_h^{-1}(i))\in S_{\vn}\quad \forall i \end{array}}\Hom(\s_l W,\s_h W)\otimes \bigoplus_i\Ext^{2}_{B\otimes B^o}\left(B,\Hom( Y_{\alpha_{\s_l^{-1}(i)}},Y_{\alpha_{\s_h^{-1}(i)}})\right)\right)^{S_N}
$$
where $(\s_l^{-1}(i) \s_h^{-1}(i))$ denotes the transposition moving the corresponding indices. Now we have $ (\s_l^{-1}(i) \s_h^{-1}(i))\in S_{\vn}$, $\forall i$ if and only if $\s_l=\s_h\s$ with $\s\in S_{\vn}$. But   $\s_l, \s_h$ belong to different left cosets of $S_{\vn}$. Thus we  can rewrite the last term as:
 $$
\left(\bigoplus_{l}\End\,\s_l W\otimes \bigoplus_i\Ext^{2}_{B\otimes B^o}(B,\End\, Y_{\alpha_{\s_l^{-1}(i)}})\right)^{S_N}=
$$
 $$
=\left(\bigoplus_{i=1}^r\End\,W\otimes \left(\Ext^{2}_{B\otimes B^o}\left(B,\End\, Y_{i}\right)\right) ^{\oplus N_i}\right)^{S_{\vn}}=
$$
$$
=\left(\bigoplus_{i=1}^r\left(\bigotimes_{j=1}^{r}\End\,W_j\right)\otimes \left(\Ext^{2}_{B\otimes B^o}\left(B,\End\, Y_{i}\right)\right)^{\oplus N_i}\right)^{S_{\vn}}.
$$

Now  as  an  $S_{N_i}$-module,  $\left(\Ext^{2}_{B\otimes B^o}\left(B,\End\, Y_{i}\right)\right)^{\oplus N_i}= \Ext^{2}_{B\otimes B^o}\left(B,\End\, Y_{i}\right)\otimes \C^{N_i}$, where $S_{N_i}$ acts only on $\C^{N_i}$ permuting the factors  and $\C^{N_i}=\C\oplus\mathfrak{h_i}$, where $\C$ is the trivial representation and $\mathfrak{h_i}$ is the reflection representation of $S_{N_i}$. So we have: $$
\bigoplus_{i=1}^r \Ext^{2}_{B\otimes B^o}\left(B,\End\,Y_{i}\right)\otimes \left(\End\,W_1\otimes\dots \otimes (\C^{N_i}\otimes \End\,W_i)\dots\otimes \End\,W_r\right)^{S_{N_1}\times\dots\times S_{N_r}}=
$$
 $$
=\bigoplus_{i=1}^r \Ext^{2}_{B\otimes B^o}\left(B,\End\,Y_{i}\right)\otimes\left(\End_{S_{N_1}}\,W_1\otimes\dots \otimes \left(\C\otimes \End\,W_i\oplus \mathfrak{h_i}\otimes \End\,W_i\right)^{S_{N_i}}\otimes \dots\otimes \End_{S_{N_r}}\,W_r\right)=
$$
 $$
=\bigoplus_{i=1}^r \Ext^{2}_{B\otimes B^o}\left(B,\End\,Y_{i}\right)\otimes \left(\End_{S_{N_i}}\,W_i\oplus \Hom_{S_{N_i}}(\mathfrak{h_i}\otimes W_i, W_i)\right)=
$$
 $$
=\bigoplus_{i=1}^r \Ext^{2}_{B\otimes B^o}\left(B,\End\,Y_{i}\right)=\C^r
$$
as $\Hom_{S_{N_i}}(\mathfrak{h_i}\otimes W_i,W_i)=0$, $\forall i$ by Lemma \ref{sym} part \emph{(i)}.

\epf

\begin{coro}\label{surj}
The map $\eta:U\longrightarrow H^2(A,\End\,M)$ is surjective.
\end{coro}
\pf\  
Let $U_0\subset U$ be the subspace of vectors $(0,c')$. 
It is sufficient to show that the restriction
of 
$\eta$ to $U_0$ is surjective. 
But this restriction is a composition of three natural maps:
$$
U_0\to H^2(B,B)\to H^2(A,A)\to H^2(A,\End\, M).
$$
Here the first map $\eta_0: U_0\to H^2(B,B)$
is induced by the deformation of $B$ along $U_0$, 
the second map $\xi: H^2(B,B)\to H^2(A,A)$ 
comes from the K\"unneth formula, 
and the third map $\psi: H^2(A,A)\to H^2(A,\End\, M)$ 
is induced by the homomorphism $A\to \End M$. 

Now, by Proposition \ref{endot}, 
the map $\psi\circ \xi$ coincides with the 
map $\psi_0: H^2(B,B)\to \bigoplus_{i=1}^r H^2(B,\End\, Y_i)$ induced 
by the homomorphism $\phi: B\to \bigoplus_{i=1}^r\End\, Y_i$.  Let  $K_0=\Ker(\psi_0)$ and  $U_0'=\eta_0^{-1}(K_0)$. We have to show that  $\mathrm{codim}\,U_0'\geq r$. But, by Proposition \ref{defo}, the representation $Y_i$ 
admits a first order deformation along $U_0'$, for any $i=1,\dots ,r$.
So on $U_0'$ we have $\Tr_{Y_i}(\Lambda)=0$, where $\Lambda$ is as in equation (\ref{centralel}) (see equation (\ref{rankone})).
Since the representation $Y_i$ correspond to the basis $\Sigma_{\lambda}$ of the root system $R_{\lambda}$ (see Theorem \ref{cbht}), 
their characters are linearly independent, hence so are these linear equations.
Thus the codimension of $U_0'$ is $\ge r$ and $\eta|_{U_0}$ is surjective, as desired.

\epf

\begin{prop}\label{h1}
$H^1(A,\End M)=0$. 
\end{prop}

\pf\ 
Arguing as in the proof of Proposition 
\ref{endot}, we get that $H^1(A, \End\, M)=\bigoplus_{i=1}^rH^1(B,\End\, Y_i)$, which is
zero. This proves the proposition. 

\epf

We have thus proved the following result. 

\begin{prop} \label{main}
If the Young diagram corresponding 
to $W_i$ is a rectangle for any $i$, $Ext_B(Y_i,Y_j)=0$ for any $i\neq j$ and $M'=W\otimes Y=W_1\otimes\dots\otimes W_r\otimes Y_1^{\otimes N_1}\otimes\dots\otimes Y_r^{\otimes N_r}$  then there exists 
a unique smooth codimension $r$ formal subscheme $S$
of the formal neighborhood of the origin in $U$
such that the representation $M=\ind _{S_{\vn}\sharp B^{\otimes N}}^{S_N\sharp B^{\otimes N}}M'$ of 
$H_{1,0,c_0}(\G_N)$ formally deforms 
to a representation of $H_{1,k,c_0+c'}(\G_N)$ 
along $S$ (i.e., abusing the language, for $(k,c')\in S$). 
Furthermore, the deformation of $M$ over $S$ is unique.
\end{prop}

\pf\ Corollary \ref{surj} and Proposition \ref{h1} show 
that our case satisfies 
all the hypotheses of Proposition \ref{defo}. Moreover, from  $H^2(A, \End\,M)=\C^r$ we deduce $\di\Ker\,\eta=\di U-r$, and the Proposition follows.

\epf

\subsection{The subscheme $S$}
 
Now we would like to find the subscheme $S$ of Proposition \ref{main}. We will do this computing some appropriate trace conditions for the deformation of the module $M$. 
To this end, we will use the definition of $H_{1,k,c}(\G_N)$ given at the end of Section \ref{wreprod} ( see Lemma \ref{pres}).
We recall that from the results of  \cite{CBH}, that we summarized in Section \ref{rank1},  the irreducible representations of $H_{1,c}(\Gamma)$ are in one-to-one correspondence with the basis $\Sigma_{\lambda}$ of the root system $R_{\lambda}$, under the assignment sending  each simple module to its dimension vector (see Theorem \ref{cbht}).

Since   $\{Y_i\}$ is a collection of non-isomorphic irreducible representations the dimension vectors attached to these modules are  distinct roots of such a basis, so in particular they are linearly independent. Thus for any complex numbers $z_1,\dots,z_r$ there exists a central element $Z$ of $\C[\Gamma]$ such that $\Tr_{Y_i}(Z)=z_i$. Fix now $Z(1)$ such that $\Tr_{Y_i}(Z(1))=1-\delta_{1i}$ and consider the element: 
$$
P_1=\underbrace{1\otimes\dots\otimes 1}_{N_1}\otimes \underbrace{Z(1)\otimes\dots\otimes Z(1)}_{N-N_1}\in\h.
$$ 
Such element commutes with $x_1$, $y_1$. Consider now the relation  \emph{(R1)} for $i=1$ and multiply it by $P_1$ on the right. The left hand side becomes $[x_1, y_1P_1]$, thus it has trace zero. To compute the trace of the right hand side operator  it is convenient to use again the decomposition of the induced module  $M$ given in the previous section (cfr. (\ref{deco}), (\ref{action})).
The trace of the left hand side reduces to the sum of the three terms:
\begin{equation}\label{trace1}
\Tr_M P_1= a\,(\di Y_1)^{N_1}\prod_{j}\di W_j
\end{equation}
\begin{equation}\label{trace2}
\frac{k}{2}\Tr_M\left(\sum_{j=2}^{N_1} \sum_{\g\in\Gamma}s_{1j}\g_1\g_j^{-1}\right)P_1=\frac{k}{2}\,a\,\Tr_{W_1\otimes Y_1^{\otimes N_1}}\left(\sum_{j=2}^{N_1} \sum_{\g\in\Gamma}s_{1j}\g_1\g_j^{-1}\right)\prod_{j\neq 1}\di W_j
\end{equation}
\begin{equation}\label{trace3}
\Tr_M\left(\sum_{\g\in\Gamma\backslash \{1\}}c_{\g}\g_1\right)P_1=a\, \left(\sum_{\g\in\Gamma\backslash \{1\}}c_{\g}\Tr_{Y_1}(\g)\right) \,(\di Y_1)^{N_1-1}\prod_j \di W_j
\end{equation}
where $a=\frac{(N-N_1)!}{N_2!\dots N_r!}$ is the number of elements in the set of representatives $\s_l\in\{\s_1,\dots,\s_n\}$ such that $\s_l^{-1}\left(\{1,2,\dots,N_1\}\right)=\{1,2,\dots,N_1\}$. 
But now we claim that:

\begin{equation}\label{trace4}
\Tr_{W_1\otimes Y_1^{\otimes N_1}}(s_{1j}\,\g_1\,{\g_j}^{-1})=\Tr_{W_1}(s_{1j})\,{\di Y_1}^{N_1-1}\qquad \forall\, j\leq N_1.
\end{equation}
To obtain (\ref{trace4}), we observe that, for $j\leq N_1$, $s_{1j}\,\g_1\,{\g_j}^{-1}$ is conjugate in $\G_{N_1}=S_{N_1}\ltimes \Gamma^{N_1}$ to $s_{1j}$ and that the character of $S_{N_1}$ on $W_1\otimes Y_1^{\otimes N_1}$ is simply the product of the characters on $W_1$ and $Y_1^{\otimes N_1}$. An easy computation gives $\Tr_{Y_1^{\otimes N^1}}s_{ij}={\di Y_1}^{N_1-1}$, hence the formula.

By  Lemma \ref{sym}, we have that for any transposition $\sigma\in S_{N_1}$, $\Tr_{W_1}(\sigma)=\frac{\di W_1}{N_1\,(N_1-1)/2} {\mathbf c}(\mu)$, where ${\mathbf c}(\mu)$ is the content of the Young diagram $\mu$ attached to $W_1$. In particular, if $\mu$ is a rectangular diagram of size $l_1\times m_1$ with $l_1m_1=N_1$, it can be easily computed that: 
$$
{\mathbf c}(\mu)=\frac{N_1\,(m_1-l_1)}{2},
$$
so we have:
\begin{equation}\label{trace5}
\Tr|_{W_1\otimes Y_1^{\otimes N_1}} (s_{1j}\g_1\,{\g_j}^{-1})=\frac{(m_1-l_1)\,\di W_1}{N_1-1} {\di Y_1}^{N_1-1}.
\end{equation}
 Substituting in (\ref{trace2}), summing up (\ref{trace1}),(\ref{trace2}),(\ref{trace3}) and simplifying we get the equation:
\begin{equation}\label{hyperplane}
\di Y_1+\frac{k}{2}|\Gamma|(m_1-l_1)+\sum_{\gamma\in
\Gamma\setminus \lbrace{1\rbrace}}c_\gamma\chi_{Y_1}(\gamma)=0
\end{equation}
that is exactly the equation for the hyperplane $\mathcal{H}_{Y_1,m_1,l_1}$. Analogously we can define $Z(i)$, $P_i$  and  obtain the equations of the  hyperplane $\mathcal{H}_{Y_i,m_i,l_i}$ for $i=2,\dots, r$. We get  exactly $r$ independent necessary linear conditions. 

This shows that $(0,c_0)+S\subset\bigcap_{i=1}^r\mathcal{H}_{Y_i,m_i,l_i}$. But since the two subschemes have the same dimension we have that $S$ is the formal neighborhood of zero in $\bigcap_{i=1}^r\mathcal{H}_{Y_i,m_i,l_i}-(0,c_0)$ and Theorem \ref{mainth}, $(i)$ is proved. Parts $(ii)$, $(iii)$ follow exactly as in \cite{EM}.

\textbf{Acknowledgments}. I am very grateful to Pavel Etingof for many comments and useful discussions.


\begin{thebibliography}{APK}




\bibitem[CBH]{CBH} W. Crawley-Boevey, M. Holland,
{\em Noncommutative deformations of Kleinian singularities},
Duke Math. J. {\bf 92} (1998), no. 3, 605--635. 

\bibitem[EG]{EG} P. Etingof, V. Ginzburg,
{\em Symplectic reflection algebras, Calogero-Moser space, and 
deformed Harish-Chandra homomorphism}, Invent. Math. {\bf 147} 
(2002), no. 2, 243--348, {\tt math.AG/0011114}.

\bibitem[EO]{EO} P. Etingof, A. Oblomkov,
{\em Quantization, orbifold cohomology, and Cherednik algebras},
preprint, {\tt math.QA/0311005}.

\bibitem[GG]{GG} W. L. Gan, V. Ginzburg 
{\em Deformed preprojective algebras and symplectic reflection algebras for wreath products},preprint, {\tt math.QA/0401038}


\bibitem[VB1]{VB1}M. Van Den Bergh
{\em A relation between Hochschild homology and cohomology for Gorenstein rings},Proc. Amer. Math. Soc. {\bf126} (1998), no. 5, 1345--1348. 

\bibitem[VB2]{VB2}M. Van Den Bergh
{\em Erratum to "`A relation between Hochschild homology and cohomology for Gorenstein rings"'},Proc. Amer. Math. Soc. {\bf130} (2000), no. 9, 2809--2810. 

\bibitem[EM]{EM} P. Etingof, S. Montarani
{\em Finite dimensional representations of symplectic reflection algebras associated to wreath products},
preprint, {\tt math.RT/0403250 v2}.




\end{thebibliography}
\end{document}